\documentclass{article}
\usepackage{graphicx} 
\usepackage[paperwidth=7in, paperheight=10in, margin=.875in]{geometry}
\usepackage[colorlinks,linkcolor=red,anchorcolor=green,citecolor=blue]{hyperref}
\usepackage{amsfonts, amssymb, amsmath, amsthm}
\usepackage{stmaryrd} 
\usepackage{graphicx} 
\usepackage{subcaption} 
\usepackage{cite}
\usepackage{cellspace} %
\usepackage{enumerate}
\usepackage[utf8]{inputenc}
\usepackage{comment}
\usepackage{longtable}
\usepackage{adjustbox}
\usepackage{float}
\usepackage{subcaption}

\usepackage[makeroom]{cancel}
\usepackage{array} 
\newcolumntype{C}[1]{>{\centering\arraybackslash}p{#1}}
\usepackage{xcolor} 
\usepackage{makecell}
\usepackage{authblk} 

\title{High-order time splitting methods for the
nonlinear Gross-Pitaevskii equation}

\author[1]{\scriptsize Roberto Ben}
\author[2]{Agustín Besteiro}
\author[3]{Diego Rial}
\affil[1]{Instituto del Desarrollo Humano, Universidad Nacional de General Sarmiento, Los Polvorines, Buenos Aires, Argentina, 
{\href{mailto:rben@campus.ungs.edu.ar}{\textcolor{blue}{rben@campus.ungs.edu.ar}}}}
\affil[2]{Centro de Matem\'{a}tica Aplicada, Instituto de Tecnolog\'{i}as Emergentes y Ciencias Aplicadas (ITECA), Universidad Nacional de San Mart\'{i}n - CONICET, Buenos Aires, Argentina,
{\href{mailto:abesteiro@unsam.edu.ar}{\textcolor{blue} {abesteiro@unsam.edu.ar}}}}
\affil[3]{Departamento de Matemática, Universidad de Buenos Aires, Instituto de Matemática Luis Santaló (IMAS) - CONICET, Argentina, 
{\href{mailto:drial@dm.uba.ar}{\textcolor{blue} {drial@dm.uba.ar}}}}

\date{\today}

\begin{document}

\maketitle

\begin{abstract}
We propose a high-order numerical methodology for computing the ground state and time evolution of the two-dimensional Gross–Pitaevskii equation with harmonic trapping potential. The ground state is obtained by combining normalized gradient flow with time-splitting schemes featuring strictly positive coefficients, which makes them suitable for both dissipative and Hamiltonian regimes and avoids the negative time steps required by classical symplectic methods. The computational cost of these methods grows quadratically with the order, in contrast to the exponential growth of symplectic alternatives. Numerical benchmarks assess ground state convergence under different initializations, long-time preservation of mass and Hamiltonian energy for varying mass constraints, and the cost–accuracy trade-off for orders $q = 2, 4, ..., 14$. 
\end{abstract}

\providecommand{\keywords}[1]
{\textbf{\textit{keywords---}} #1}

\providecommand{\AMS}[1]
{\textbf{\textit{AMS code---}} #1}

\section{Introduction}
The Gross-Pitaevskii equation (GPE) models the evolution of the macroscopic wave function of a Bose-Einstein condensate (BEC) under the influence of an external potential and interactions between particles in the condensate \cite{pethick2008bose, pitaevskii1961vortex, gross1961structure}. Its general formulation is as follows:
\begin{align}\label{eq:Evol_Gross-Pitaevskii}
\begin{split}
    i \partial_t \psi(\mathbf{x},t) & =  -\frac{1}{2} \nabla^{2}\psi(\mathbf{x},t)
    + \left( V(\mathbf{x},t) + \beta |\psi(\mathbf{x},t)|^{2}\right) \psi(\mathbf{x},t), \, \mathbf{x} \in \Omega, t \in [0,T], \\
    \psi(\mathbf{x},t_{0})  & = \psi_{0}(\mathbf{x}),
\end{split}
\end{align}
where $\psi(\mathbf{x},t)$ is a complex-valued function, $\beta > 0$, $\Omega \subseteq \mathbb{R}^d$, $d \in \left\{1,2,3\right\}$ and $T>0$.
We aim to numerically solve this initial value problem and to find the ground state of the GPE. We focus on the two-dimensional GPE with an isotropic harmonic
potential
\[
V(\mathbf{x},t)=\frac{1}{2}\gamma^2 (x_1^2+x_2^2),
\]
$\gamma >0$, which is a standard model in the description of a magnetically trapped BEC in the quasi-two dimensional regime \cite{kevrekidis2015defocusing}.
From a mathematical and numerical viewpoint, this choice of the potential is particularly convenient, since the linear part of the Hamiltonian admits an explicit spectral decomposition in terms of Hermite functions. As a consequence, the associated linear flow can be
computed exactly within a truncated Hermite basis, making spectral time-splitting
methods especially effective.

Although the splitting schemes considered in this work are not conservative by
construction, our numerical experiments indicate that in one spatial dimension ($d=1$) the $L^{2}$-norm and the Hamiltonian are preserved up to machine precision. We do not include the one-dimensional case here in order to avoid additional material that is not essential to the main objectives of this paper. In contrast, in two spatial dimensions, the Fourier--Hermite transforms involve genuinely two-dimensional matrix operations, whose truncation and numerical implementation introduce non-negligible error propagation. This makes the two-dimensional setting particularly suitable for assessing the stability, accuracy, and conservation properties of high-order time-splitting methods, and motivates our focus on the case $d=2$.
Finally, we note that the proposed methodology is not restricted to two spatial dimensions and can be extended to the three-dimensional GPE.

High-order splitting techniques were developed based on Lie-Trotter schemes for the numerical approximation of solutions \cite{Borgna2015,DeLeo2015}. Recently, these splitting methods were used to study well-posedness of parabolic equations  \cite{besteiro1,Besteiro2018}.
In the historical continuum of quantum physics, the GPE stands as a quintessential mathematical model. For several decades, this equation has been important as the theoretical framework for describing Bose-Einstein condensed quantum systems. Originally formulated to elucidate the behavior of dilute atomic gases at ultracold temperatures, this equation has found extensive applications in various scientific disciplines, including atomic and molecular studies, solid-state physics, nuclear physics, nonlinear optics and chemistry \cite{lieb2001rigorous, pelinovsky2011localization, pitaevskii1961vortex, borgna2018optical, assanto2012nematicons,ben2015localized,ben2017properties, kevrekidis2008emergent, carretero2024nonlinear}. 
Over the years, numerous works have successfully developed numerical solutions of the GPE employing diverse computational methods \cite{bao2009generalized, antoine2013computational, fu2022arbitrary}.
This paper introduces a numerical methodology for solving the two-dimensional GPE. 
The ground state is computed by combining splitting methods, gradient descent, and $L^2$-normalization at each iteration. The time evolution is then studied by initializing the GPE with the computed ground state and applying the same splitting-based integrators. A comprehensive set of numerical benchmarks is carried 
out, addressing ground state convergence under different initializations, long-time 
preservation of mass and Hamiltonian energy, and the cost--accuracy trade-off for 
method orders $q = 2, 4, \ldots, 14$.

The paper is organized as follows. In Section \ref{sec:High-order_TSM} we briefly outline the positive time-stepping high-order splitting methods proposed in~\cite{DeLeo2015}, emphasizing their stability, applicability to irreversible problems, computational efficiency, and parallelization capabilities. In Section \ref{sec:Groundstate} we present the splitting approach for ground-state computation and analyze the associated linear and nonlinear subproblems. In Section \ref{sec:gpedynamics} we apply the splitting strategy to the real-time dynamics of the GPE, again describing the corresponding linear and nonlinear substeps. Finally, in Section \ref{sec:numericalapprox} we report numerical results: we study the ground-state computation under different initializations and their convergence behavior. We then investigate the time evolution of the computed ground state, with particular emphasis on the numerical conservation of mass and energy. The resulting relative errors, together with the corresponding computational times, are consistent with those reported in \cite{fu2022arbitrary}, where splitting methods up to order six are considered for the free case $V=0$. The present work goes beyond that setting by treating methods of order up to fourteen in the presence of a harmonic trapping potential.

\section{High-order time splitting methods for the GPE}\label{sec:High-order_TSM}


The initial value problem (\ref{eq:Evol_Gross-Pitaevskii}) is a particular case of the general problem:

\begin{equation*}\label{eq:initial_value_problems}
    \partial_t \psi = A_0 \psi + A_1(\psi), \quad \psi(0) = \psi_0,
\end{equation*}
where \(A_0\) is a linear closed operator with dense domain \(D(A_0) \subset \mathbf{H}\), \(\mathbf{H}\) is a Hilbert space, and \(A_1: \mathbf{H} \rightarrow \mathbf{H}\) is a smooth mapping with \(A_1(0) = 0\). When the partial differential equations \(u_t = A_0 u\) and \(u_t = A_1(u)\) can be solved independently, approximate solutions to the problem can be obtained by alternately applying the flows associated with each problem. In De Leo et al. \cite{DeLeo2015}, the consistency, stability, and convergence of a collection of high-order methods for irreversible equations are presented. These methods are based on the ideas introduced in the Lie-Trotter and Strang numerical integration methods, defined by
\begin{align*}
    \Phi_{\text{Lie}}(\tau,u_0) &= \phi_1(\tau,\phi_0(\tau,u_0)), \\
    \Phi_{\text{Strang}}(\tau,u_0) &= \phi_0\left(\frac{\tau}{2},\phi_1 \left(\tau,\phi_0\left(\frac{\tau}{2},u_0\right)\right)\right).
\end{align*}
where \(\tau\) is the time step of numerical integration, and \(\phi_0\) and \(\phi_1\) are the flows associated with the partial differential equations \(\partial_t u = A_0u\) and \(\partial_t u = A_1(u)\), respectively. It is a well-known fact that $\Phi_{\text{Lie}}$ has order one and $\Phi_{\text{Strang}}$ has order two. Methods of order $q=3,4,2n$ were presented by Ruth \cite{ruth1983canonical}, Neri \cite{neri1987lie}, and Yoshida \cite{yoshida1990construction}.

In this century many authors presented rigorous results on the convergence of symplectic methods applied to Hamiltonian systems of infinite dimensions, e.g., Besse et al.~\cite{besse2002order}, Descombes \& Thalhammer \cite{descombes2010exact, descombes2013lie} for the cubic nonlinear Schr\"odinger equation, and Lubich~\cite{lubich2008splitting} and Gauckler~\cite{gauckler2011convergence}
for the GPE. For orders $q>2$, symplectic methods necessarily involve negative time steps, which inhibits their application to irreversible evolution problems.
This limitation is particularly relevant in the computation of the ground state via normalized gradient flows, which generate only a forward-time semigroup. For this reason, in the following we focus on the class of high-order time-splitting methods introduced in~\cite{DeLeo2015}, whose coefficients are strictly positive, and we specify which of them are employed in Section~\ref{sec:Groundstate} to approximate the solution of Eq.~\eqref{eq:Evol_Gross-Pitaevskii}.

Let $\Phi^{+}_1(\tau)= \phi_1(\tau) \circ \phi_0(\tau)$ and $\Phi^{-}_1(\tau)= \phi_0(\tau) \circ \phi_1(\tau)$ and, for $m>1$, $\Phi^{\pm}_m(\tau) = \Phi^{\pm}_1(\tau) \circ \Phi^{\pm}_{m-1}(\tau)$. In \cite{DeLeo2015} the authors describe two types of time-splitting methods of arbitrary even order $q=2n$. In the following, we focus on the symmetric method defined by:
\begin{equation}\label{eq:symmetric_time_splitting_method}
    \Phi(\tau) = \sum_{m=1}^{s} \gamma_m \left(\Phi^{+}_{m}(\tau/m) + \Phi^{-}_{m}(\tau/m)\right),
\end{equation}
where the coefficient $\gamma_m$ ($m = 1,\dots, s$) satisfy the conditions:
\begin{equation}\label{eq:conditions_sym_TSM}
\begin{aligned}
    \frac{1}{2} &= \gamma_1 + \gamma_2 + \cdots + \gamma_s, \\
    0 &= \gamma_1 + 2^{-2k}\gamma_2 + \cdots + s^{-2k}\gamma_s, \quad 1 \leq k \leq n-1.
\end{aligned}
\end{equation}
The total number of basic steps is given by $S = 4 \sum_{\gamma_m \neq 0} m$. For example, for $s=n$ the system (\ref{eq:conditions_sym_TSM}) has a unique solution with $\gamma_m \neq 0$ for $1\leq m \leq s$, and the total number of steps is given by
\begin{equation}\label{eq:steps_order}
S = q(\frac{q}{2}+1), 
\end{equation}
which grows quadratically with order.  In contrast, the symplectic methods presented in Yoshida \cite{yoshida1990construction} and Castella et al. \cite{castella2009splitting} exhibit exponential growth, highlighting their distinct computational complexity. Moreover, the ability to compute $\Phi_m^{\pm}$ in parallel significantly reduces the total computation time when using multiple processors.

Next, we show how these time-splitting methods can be adapted and extended to compute the ground state and to compute the subsequent dynamics of the GPE (\ref{eq:Evol_Gross-Pitaevskii}). In particular, we incorporate an $L^{2}$-normalization step within each iteration of the gradient flow used for ground state computation.
This modification is essential, since the underlying descent dynamics does not preserve the mass constraint.
To this end, we briefly recall the formulation of the GPE and
introduce the specific splitting-based schemes employed in this work.

\section{Groundstate approximation}
\label{sec:Groundstate}
It is well-known that the GPE \cite{kevrekidis2015defocusing}:
\begin{equation}\label{eq:2D_Evol_Gross-Pitaevskii}
i \partial_t \psi(\mathbf{x},t) = -\frac{1}{2} \nabla^{2}\psi(\mathbf{x},t)
+ \left( V(\mathbf{x}) + \beta |\psi(\mathbf{x},t)|^{2}\right) \psi(\mathbf{x},t), \, \mathbf{x} \in \Omega, t \in [0,T],
\end{equation}
where $\psi(\mathbf{x},t)$ is a complex-valued function, $\beta > 0$, $\Omega \subseteq \mathbb{R}^d$, $d \in \left\{1,2,3\right\}$ and $T>0$, is formally equivalent to the Hamiltonian system:
\begin{align}\label{eq:Hamiltonian_Evol_Gross-Pitaevskii}
\partial_t \psi & = -i\frac{\partial \mathcal{H}}{\partial \psi^*}, \ \text{with} \ \ \ 
\mathcal{H}(\psi) = \frac{1}{2} \int_{\Omega}\left\vert \nabla\psi \right\vert ^{2} + \int_{\Omega} V \left\vert \psi \right\vert ^{2} + \frac{1}{2} \int_{\Omega} \beta \left\vert \psi \right\vert ^{4}.
\end{align}
In what follows $\Omega = \mathbb{R}^d$. The $L^{2}(\Omega)$ norm of $\psi$, $\Vert\psi(t)\Vert_{L^2} := \left(\int_\Omega |\psi(\mathbf{x},t)|^2\,d\mathbf{x}\right)^{1/2}$, and the Hamiltonian $\mathcal{H}(\psi(t))$, are conserved quantities; i.e., $\frac{d}{dt}\Vert\psi(t)\Vert_{L^2} = \frac{d}{dt}\mathcal{H}(\psi(t)) = 0$.

\label{subsec:groundstateapp}
We consider stationary solutions of Eq.~\eqref{eq:2D_Evol_Gross-Pitaevskii} of the form $\psi(\mathbf{x},t) = e^{-i\mu t}u(\mathbf{x})$, where $\mu \in \mathbb{R}$ is the so-called chemical potential associated with the temporal frequency of the solution. These solutions satisfy the stationary equation:
\begin{equation}\label{eq:stationary_GPE}
-\frac{1}{2} \nabla^{2} u(\mathbf{x})
+ \left( V (\mathbf{x}) + \beta |u(\mathbf{x})|^{2}\right) u(\mathbf{x}) = \mu u(\mathbf{x}).
\end{equation}
A special stationary solution is the groundstate, which minimizes $\mathcal{H}$ and is a solution of the Lagrangian problem:
\begin{equation*}\label{eq:H_minimizer}
\left\{
\begin{array}{cl}
\min\limits_{\psi} \mathcal{H}(\psi) \\
\Vert \psi \Vert_{L^2} = \text{constant}
\end{array}
\right. .
\end{equation*}
We use the gradient descent method combined with high-order time-splitting methods to numerically compute the ground state. The gradient descent method is defined by the equation:
\begin{equation}\label{eq:gradient_descent_method} 
    \partial_t \psi(\mathbf{x},t) = \frac{1}{2} \nabla^{2}\psi(\mathbf{x},t) - V(\mathbf{x}) \psi(\mathbf{x},t) - \beta \vert \psi(\mathbf{x},t) \vert^2 \psi(\mathbf{x},t), 
\end{equation}
on the sphere $\left\{u \in L^2(\mathbb{R}^d): \Vert u \Vert_{L^2}=c \right\}$ (see \cite{bao2004computing} for the standard construction of the equation). 
As mentioned above, we consider $d=2$ and the potential $V(x,y) = \frac{1}{2} \, \gamma^{2} \left(x^{2} + y^{2}\right)$. Observe that, due to the first term on the right-hand side of Eq.~\eqref{eq:gradient_descent_method}, this is a nonlinear diffusion equation, for which symplectic methods with negative time steps cannot be applied. Therefore, the high-order methods described above in Eq.~\eqref{eq:symmetric_time_splitting_method} are suitable. Specifically, we apply these methods by applying a time-splitting to Eq.~\eqref{eq:gradient_descent_method} considering the partial flows associated with the linear and nonlinear terms. The linear problem \eqref{eq:groundstate_linear} is solved using a spectral method, while the nonlinear problem (\ref{eq:groundstate_nonlinear}) is addressed by evolving its exact solution. Namely, we will solve separately the following two problems:
\begin{align}
\label{eq:groundstate_linear} 
\partial_t \psi(x,y,t) & = \frac{1}{2} \nabla^{2}\psi(x,y,t)
- \frac{1}{2} \, \gamma^{2} \left(x^{2} + y^{2}\right) \psi(x,y,t),  \\
\label{eq:groundstate_nonlinear}  
\partial_t \psi(x,y,t) & = -\beta \, |\psi(x,y,t)|^{2} \, \psi(x,y,t).
\end{align}
To preserve the $L^2$-norm, we adapt these methods by normalizing the solution
after each step of the high-order time-splitting procedure.

\subsection{Approximation of the linear and nonlinear flows}\label{subsec:Linear_term_approximation}

To solve Eq.~\eqref{eq:groundstate_linear}, we apply a spectral method, leveraging the fact that the Hermite functions, $$h_{k,l}(x,y) = \frac{1}{\sqrt{2^{k} \, k!2^{l} \, l!}} \, \left(\frac{\gamma}{\pi}\right)^{1/2}\, e^{-\frac{1}{2} \, \gamma \, \left(x^{2}+y^{2}\right)} \, H_{k}\left(\gamma \, x\right)H_{l}\left(\gamma \, y\right), \, \, k,l \in \mathbb{N}\cup\left\{0\right\},$$ where $H_{n}\left(z\right)$ is the $n$-th Hermite polynomial, form an orthogonal basis of eigenfunctions for the linear operator ${\cal{L}}(u) :=\frac{1}{2} \nabla^{2}u(x,y,t)
- \frac{1}{2} \, \gamma^{2} \left(x^{2} + y^{2}\right) u(x,y,t)$. Namely,
\begin{align}\label{eq:2D_Hermite_eigenfunctions}
\frac{1}{2} \nabla^{2}h_{k,l}(x,y)- \frac{1}{2} \, \gamma^{2} \left(x^{2} + y^{2}\right) h_{k,l}(x,y) = -\mu_{k,l} \, h_{k,l}(x,y), \, \, k,l \in \mathbb{N}\cup\left\{0\right\},
\end{align}
where $\mu_{k,l}=\gamma (k+l+1)$ is the eigenvalue associated to $h_{k,l}$.
Thus, the solution to the linear problem
\begin{align}\label{eq:linear_problem} 
\begin{aligned}
\partial_t \psi  & = \frac{1}{2} \nabla^{2}\psi
- \frac{1}{2} \, \gamma^{2} \left(x^{2} + y^{2}\right) \psi, \\
\psi(x, y,t_0) & = \psi_0(x, y),
\end{aligned}.
\end{align}
can be computed as $\psi(x,y,t) = \displaystyle{\sum_{k = 0}^{\infty}} \, \displaystyle{\sum_{l = 0}^{\infty}} e^{- \mu_{k,l} (t-t_{0})} \, a_{k,l} \, h_{k,l}(x,y)$, where
$a_{k,l} = \langle h_{k,l}|\psi_{0}\rangle$.
In the examples presented in this work, we numerically approximate the solution by projecting it onto the finite-dimensional vector space spanned by the basis ${\mathbb{H}}_{M,2} := \left\{ h_{k,l}(x,y) : 0 \leq k, l \leq M \right\}$, where $M$ denotes the truncation index:
\begin{equation*}\label{eq:linear_term_num_solution}
    \tilde{\psi}(x,y,t) = \displaystyle{\sum_{k = 0}^{M}} \, \displaystyle{\sum_{l = 0}^{M}} e^{- \mu_{k,l} (t-t_{0})} \, a_{k,l} \, h_{k,l}(x,y).
\end{equation*}
To compute $a_{k,l}$ we consider the $M+1$ zeros of the Hermite polynomial $H_{M+1}$: $\hat{z}_0, \hat{z}_1, \dots, \hat{z}_{M}$, and  the respective weights $$\hat{w}_r = \frac{\sqrt{\pi}2^{M} (M+1)! e^{\hat{z}_r^2}}{(M+1)^2 H_M^2(\hat{z}_r)}, \text{ for } 0\leq r \leq M,$$  corresponding to the Gauss-Hermite quadrature rule. These zeros and weights scale as $z_r = \gamma^{-1/2}\hat{z}_r$ and $w_r = \gamma^{-1/2}\hat{w}_r$, for $0\leq r \leq M$. The coefficients $a_{k,l}$ are computed using Gaussian quadrature:
\begin{equation}\label{eq:gaussian_quadrature}
    a_{k,l}=\langle h_{k,l}, \psi_{0} \rangle = 
\int_{-\infty}^{\infty}\int_{-\infty}^{\infty} h_{k,l}(x,y) \, \psi_0(x,y) \, dx dy
= \sum_{r = 0}^{M} \sum_{s = 0}^{M} h_{k}(z_{r}) \, h_{l}(z_{s}) \, \psi_{r,s} \, w_{r} \, w_{s},
\end{equation}
where $$h_n(z) = \frac{1}{\sqrt{2^n n!}}\left(\frac{\gamma}{\pi}\right)^{\frac{1}{4}} e^{-\frac{1}{2}\gamma z^2} H_n(\gamma z),$$ 
and $\psi_{r,s}=\psi_{0}(z_r,z_s)$.
The coefficients $a_{k,l}$ can be compactly expressed in matrix form as $a_{k,l}=\left(G\cdot \Psi_0 \cdot G^{*}\right)_{k,l}$, where $\Psi_{0_{r,s}} = \psi_0(z_r,z_s)$, $G$ is the matrix given by
\begin{align}\label{eq:G_matrix}
G = 
\left(
\begin{array}{ccc}
h_{0}(z_{0}) \, w_{0} & \cdots & h_{0}(z_{M}) \, w_{M} \\
\vdots & & \vdots \\
h_{M}(z_{0}) \, w_{0} & \cdots & h_{M}(z_{M}) \, w_{M}
\end{array}
\right),
\end{align}
and $G^*$ denotes the adjoint (i.e., conjugate transpose) of $G$.

On the other hand, the exact solution to the nonlinear partial problem 
\begin{align}
\label{eq:nonlinear_problem}  
\begin{aligned}
\partial_t \psi(x,y,t) & = -\beta \, |\psi(x,y,t)|^{2} \, \psi(x,y,t), \\
\psi(x, y,t_0) & = \psi_0(x, y),
\end{aligned}
\end{align}
is given by
\begin{align*}
\psi(x,y,t) = \frac{\psi_{0}(x,y)}{\sqrt{1 + 2 \, \beta \, (t - t_{0}) \, 
|\psi_{0}(x,y)|^{2}}}.
\end{align*}

\subsection{High-order splitting scheme for the normalized gradient flow}
\noindent
The gradient flow defined by Eq.~\eqref{eq:gradient_descent_method} does not preserve the $L^2$-norm of the solution. However, the ground state is defined as a minimizer of the Hamiltonian $\mathcal{H}$ under the mass constraint $\|\psi\|_{L^2}=c$. Consequently, the normalization step is not merely a numerical stabilization technique, but an essential component of the method, enforcing the constraint at the discrete level.
From a variational perspective, this procedure can be interpreted as a
projection of the gradient flow onto the $L^2$-sphere of radius $c$, yielding a
\emph{normalized gradient flow} in the sense of \cite{bao2004computing, BaoCai2013}. In
particular, without this normalization, the descending flow would not be
consistent with the constrained minimization problem defining the ground
state.

Within this framework, we track the evolution of the normalized gradient flow
obtained by combining high-order time-splitting schemes for the linear and
nonlinear partial flows with $L^2$ normalization at each iteration. The
numerical experiments (see Figs.~\ref{fig:example1_Hamiltonian} and \ref{fig_example2_Hamiltonian_faster}) show that this descending method converges towards a
minimizer of the Hamiltonian under the prescribed mass constraint.

\subsection{Chemical potential and phase--induced periodicity}\label{subsec:mu_period}

As mentioned above, let $u(\mathbf{x})$ be a stationary profile such that $
\psi(\mathbf{x},t)=e^{-i\mu t}u(\mathbf{x})$,
$ \|u\|_{L^2}=c,$
so that $u$ satisfies the stationary Gross--Pitaevskii Eq.~\eqref{eq:stationary_GPE}.
Multiplying 
both sides of that equation by $u(\mathbf{x})^{*}$ and integrating over $\Omega$ we obtain
\begin{equation}\label{eq:mu_id_step1}
\int_{\Omega}\left(-\frac12(\nabla^{2}u)\,u^{*}
+V|u|^{2}+\beta|u|^{4}\right)\,d\mathbf{x}
=\mu\int_{\Omega}|u|^{2}\,d\mathbf{x}
=\mu\,c^{2}.
\end{equation}
The Laplacian term can be handled by integration by parts (assuming the boundary term vanishes)
and Eq.~\eqref{eq:mu_id_step1} can be rewritten to define the chemical potential $\mu_{ch}$ as
\begin{equation}\label{eq:mu_formula_final}
\mu_{\mathrm{ch}}
:=\frac{1}{c^{2}}\left(
\frac12\int_{\Omega}|\nabla u|^{2}\,d\mathbf{x}
+\int_{\Omega}V(\mathbf{x})|u|^{2}\,d\mathbf{x}
+\beta\int_{\Omega}|u|^{4}\,d\mathbf{x}
\right).
\end{equation}
Note that the right-hand side involves the same integrals as $\mathcal{H}(u)$ (see Eq.~\eqref{eq:Hamiltonian_Evol_Gross-Pitaevskii}), and that the chemical potential $\mu_{\mathrm{ch}}$ is determined through the mass constraint $\Vert u\Vert_{L^2}=c$.

The temporal frequency determined by the chemical potential induces periodicity in the exact evolution of a stationary state. This can be seen from
\begin{align*}\label{eq:L2_distance_phase}
\|\psi(\cdot,t)-u\|_{L^2}
=\big|e^{-i\mu_{\mathrm{ch}} t}-1\big|\,\|u\|_{L^2}
=2\|u\|_{L^2(\Omega)}\left|\sin\left(\frac{\mu_{\mathrm{ch}} t}{2}\right)\right|,
\end{align*}
which yields the relation
\begin{equation}\label{eq:period_mu}
T_{\mu}:=\frac{2\pi}{\mu_{\mathrm{ch}}}.
\end{equation}

In the experiments conducted in Section~\ref{subsec:Invariant_preservation}, 
we compute $\mu_{\mathrm{ch}}$ and $T_{\mu}$, and numerically verify that the theoretical value of $T_{\mu}$ agrees with the observed period of the numerical solution.


\section{GPE dynamics}
\label{sec:gpedynamics}
In this section we study the numerical evolution of the GPE (\ref{eq:2D_Evol_Gross-Pitaevskii})
\begin{align}\label{eq:gpe_dynamics}
i \partial_t \psi(\mathbf{x},t) & = -\frac{1}{2} \nabla^{2}\psi(\mathbf{x},t)
+ \left( V(\mathbf{x},t) + \beta |\psi(\mathbf{x},t)|^{2}\right) \psi(\mathbf{x},t), \, \mathbf{x} \in \Omega, t \in [0,T],
\\
 \psi(\mathbf{x},t_{0}) & = \psi_{0}(\mathbf{x}).  \notag 
\end{align}
As we did in the last section, we can apply the time-splitting methods described in Section \ref{sec:High-order_TSM}, separating Eq.~\eqref{eq:gpe_dynamics} into the linear and nonlinear terms. Once more, the linear problem (\ref{eq:gpe_dynamics_linear}) is solved using a spectral method, while the nonlinear problem (\ref{eq:gpe_dynamics_nonlinear}) is addressed by evolving its exact solution. Note that the partial problems:
\begin{align}
\label{eq:gpe_dynamics_linear} 
i \partial_t \psi(x,y,t) & = - \frac{1}{2} \nabla^{2}\psi(x,y,t)
+ \frac{1}{2} \, \gamma^{2} \left(x^{2} + y^{2}\right) \psi(x,y,t), \\
\label{eq:gpe_dynamics_nonlinear}  
i \partial_t \psi(x,y,t) & = -\beta \, |\psi(x,y,t)|^{2} \, \psi(x,y,t), 
\end{align}
are not the same as those stated in Eq.~\eqref{eq:groundstate_linear} and Eq.~\eqref{eq:groundstate_nonlinear}.
The key difference between this system and the one considered in the previous section lies in the presence of the imaginary unit multiplying the time derivative. While the gradient flow equation defines a dissipative evolution, the GPE corresponds to a Hamiltonian (unitary) dynamics, leading to fundamentally different qualitative properties.

\subsection{Linear and nonlinear evolution subproblems}
As shown in Section \ref{subsec:Linear_term_approximation}, we can make use of the Hermite functions to find the solution of the linear equation
\begin{align*}
i \, \partial_t \psi(x,y,t) & = -\frac{1}{2} \nabla^{2}\psi(x,y,t)
+ \frac{1}{2} \, \gamma^{2} \left(x^{2} + y^{2}\right) \psi(x,y,t) \\[4pt]
\psi(x,y,t_{0}) & = \ \psi_{0}(x,y)
.
\end{align*}
The solution to this linear problem can be expressed as the complex function $$\psi(x,y,t) = \displaystyle{\sum_{k = 0}^{\infty}} \, \displaystyle{\sum_{l = 0}^{\infty}} e^{-i \mu_{k,l} (t-t_{0})} \, a_{k,l} \, h_{k,l}(x,y),$$ where
$a_{k,l} = \langle h_{k,l},\psi_{0}\rangle$.
Projecting it onto the finite-dimensional vector space spanned by the basis ${\mathbb{H}}_{M,2} := \left\{ h_{k,l}(x,y) : 0 \leq k, l \leq M \right\}$, where $M$ determines the truncation index, we can numerically approximate the evolution function of the linear term:
\begin{equation*}\label{eq:evolution_linear_term_num_solution}
    \tilde{\psi}(x,y,t) = \displaystyle{\sum_{k = 0}^{M}} \, \displaystyle{\sum_{l = 0}^{M}} e^{- i\mu_{k,l} (t-t_{0})} \, a_{k,l} \, h_{k,l}(x,y),
\end{equation*}
with
$a_{k,l} = \langle h_{k,l},\psi_{0}\rangle$, $0\leq k, l \leq M$. These coefficients can be computed using $a_{k,l}=\left(G \cdot \Psi_0 \cdot G^{*}\right)_{k,l}$, as shown in Eqs.~\eqref{eq:gaussian_quadrature} and \eqref{eq:G_matrix} in Section \ref{subsec:Linear_term_approximation}.

Furthermore, for $\psi_0$ in the vector space spanned by ${\mathbb{H}}_{M,2}$, we theoretically verify the conservation of the norm:
\begin{align*}
\langle\tilde{\psi},\tilde{\psi}\rangle & = 
\int_{-\infty}^{+\infty}\int_{-\infty}^{+\infty}
\overline{\tilde{\psi}(x,y,t)} \, \tilde{\psi}(x,y,t) \, dx \, dy \\
& = \sum_{k,l,p,q = 0}^{M}
e^{-i (\mu_{k} - \mu_{p}) \, (t-t_0)} \,
e^{-i (\mu_{l} - \mu_{q}) \, (t-t_0)} \,
\bar{a}_{k,l} \, a_{p,q} \, \langle h_{k,l} , h_{p,q} \rangle \\ 
& = \sum_{k,l = 0}^{M} |a_{kl}|^{2}
= \langle \psi_{0},\psi_{0}\rangle.
\end{align*}

For the nonlinear subproblem, we analyze the dynamics of non-linear evolution. The governing equation is given by
\begin{align*}
i \, \partial_t \psi(x,y,t) & =  - \beta \, |\psi(x,y,t)|^{2} \, \psi(x,y,t) \\[2pt]
\psi(x,y,t_{0}) & = \psi_{0}(x,y)
\end{align*}
This equation admits an exact solution since $|\psi(x,y,t)|^2$ remains constant in time. Indeed, for fixed $(x,y)$,
\[
\frac{d}{dt}|\psi|^2
= \partial_t\psi\,\overline{\psi} + \psi\,\partial_t\overline{\psi}
= -i\beta|\psi|^2\psi\overline{\psi}
+ i\beta|\psi|^2\psi\overline{\psi}
= 0.
\]
Hence,
\[
|\psi(x,y,t)|^2 = |\psi(x,y,t_0)|^2 \quad \text{for all } t,
\]
and the solution evolves through a pure phase rotation,
\[
\psi(x,y,t) = e^{-i \beta \, (t - t_{0}) \, \rho(x,y)} \, \psi_{0}(x,y),
\]
where $\rho(x,y) = |\psi_{0}(x,y)|^{2}$.

\section{Numerical results}
\label{sec:numericalapprox}

We present a set of benchmark numerical experiments aimed at providing a systematic and reproducible assessment of the proposed high-order time-splitting methods. The benchmarks address ground state convergence, sensitivity to the initial seed, long-time preservation of mass and Hamiltonian energy, and the cost–accuracy trade-off associated with different splitting orders, time steps, and CPU times.

Throughout this section, we fix $\beta=2$ and $\gamma=1$ without loss
of generality, as these parameters can be effectively absorbed by a
suitable rescaling of time and wavefunction. Unless stated otherwise,
we work with Hermite functions in ${\mathbb{H}}_{16,2}$ evaluated on
the square grid $\mathcal{T}:=\left\{(z_k,z_l)\right\}_{0\leq k,l\leq 16}$,
where $z_0,\dots,z_{16}$ are the roots of the Hermite polynomial $H_{17}$.

\subsection{Benchmark I: Ground state convergence}
\subsubsection{Sensitivity to the initial seed I}\label{Subsubsec:Example1}
\begin{figure}[h!]
\begin{subfigure}[h]{0.49\linewidth}
\includegraphics[width=\linewidth]{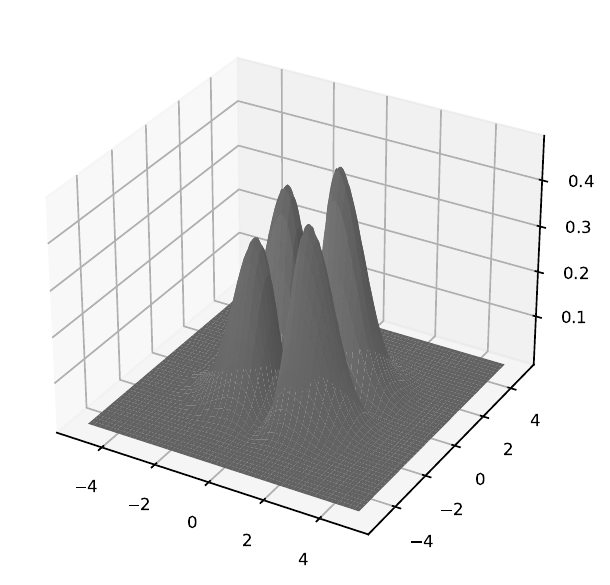}
\caption{Initial seed employed in the gradient descent method.}\label{fig:GS_initial_seed}
\end{subfigure}
\hfill
\begin{subfigure}[h]{0.49\linewidth}
\includegraphics[width=\linewidth]{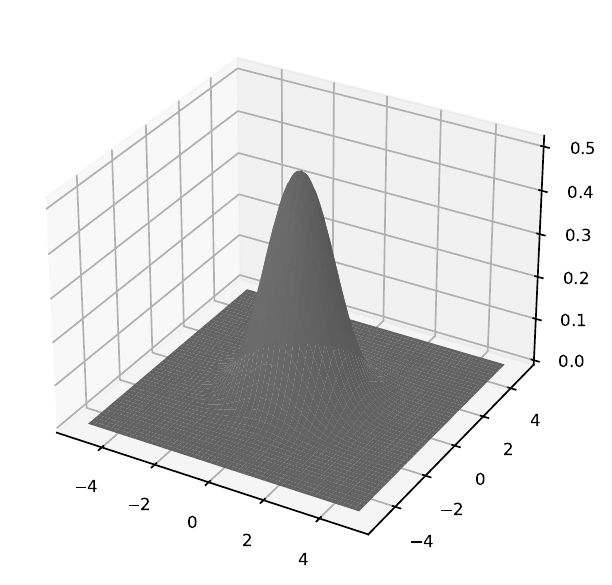}
\caption{Numerical approximation of the ground state.}
\label{fig:Groundstate}
\end{subfigure}%
\caption{Ground state computed initialized with the seed $\Psi_{0}$ with a time-splitting gradient descent method.}
\end{figure}

First, we apply the gradient descent method described in Section \ref{subsec:groundstateapp}.   The time step is set to \(\tau = 0.01\), and we start, for instance, from an initial seed (Fig.~\ref{fig:GS_initial_seed}) defined by \(\Psi_0 = 0.01 \tilde{h}_{0,0} +0.1 \tilde{h}_{0,1} + 0.1 \tilde{h}_{1,0} + \tilde{h}_{1,1}\), where $\tilde{h}_{k,l}$ represents the restriction $h_{k,l}\big|_{\mathcal{T}}$ with \(k,l \in \{0,1\}\).
Using a gradient descent method of order 4, we found the numerical ground state shown in Fig.~\ref{fig:Groundstate}.

\begin{figure}[h!]
\centering
\begin{subfigure}[h]{0.48\linewidth}
\includegraphics[width=\linewidth]{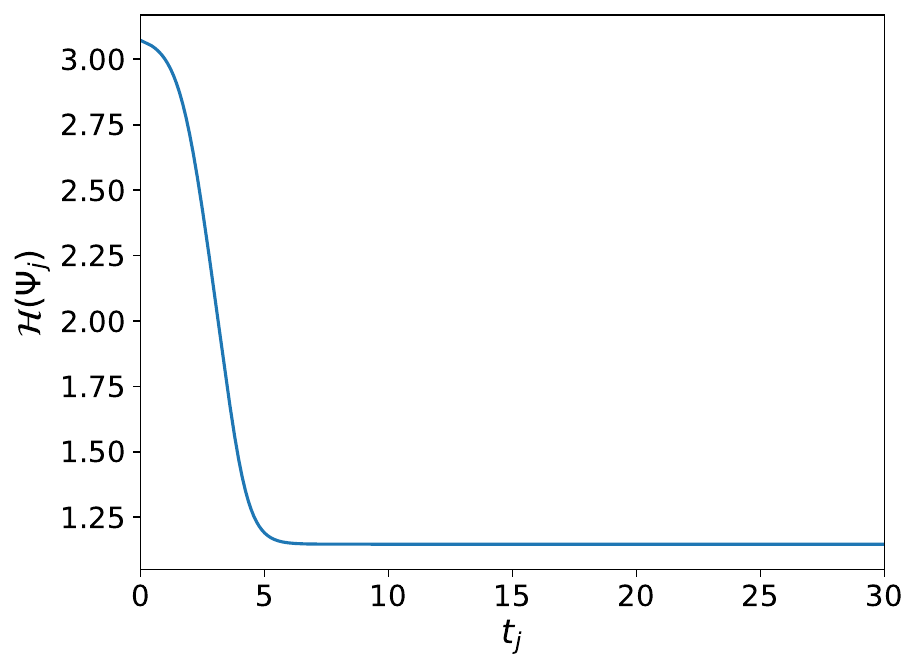}
\caption{Hamiltonian starting from initial seed of Fig.~\ref{fig:GS_initial_seed}.}\label{fig:example1_Hamiltonian}
\end{subfigure}
\hfill
\begin{subfigure}[h]{0.49\linewidth}
\includegraphics[width=\linewidth]%
{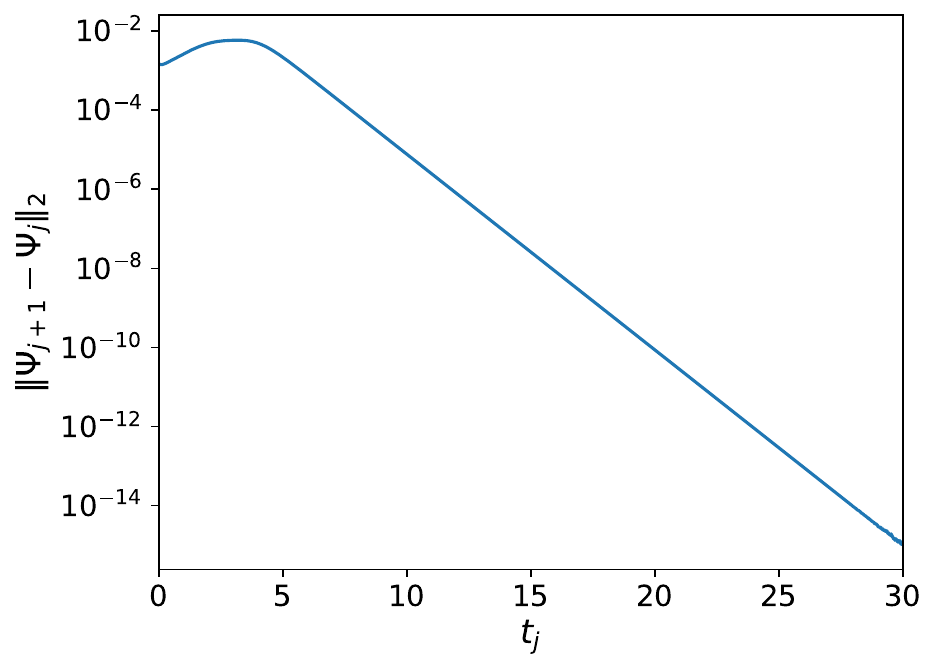}

\caption{Error between consecutive approximations of the gradient descent method.}\label{fig:example1_step_error}
\end{subfigure}
\caption{Evolution of Hamiltonian and convergence of the gradient-descent iterations initialized with the seed $\Psi_{0}$.}

\end{figure}

The initial seed chosen is close to the Hermite function \(h_{1,1}\), which is an eigenfunction of the linear part of the Hamiltonian (cf. Eq.~\eqref{eq:2D_Hermite_eigenfunctions}). Consequently, in the early time steps of the gradient descent method, we observe a slow decay of the Hamiltonian, as depicted in Fig.~\ref{fig:example1_Hamiltonian}. Subsequently, there is a rapid decay towards the Hamiltonian minimum, which is computed as \(\mathcal{H}_{\text{min}} =  1.14672998984857\), for a final simulation time $T_f = 30$ seconds.  

Let $\Psi_j$ denote the  numerical approximation after $j$ iterations of the evolution operator using the time-splitting method defined in Section $\ref{sec:Groundstate}$.
Fig.~\ref{fig:example1_step_error} displays the evolution of the distance between consecutive approximations, $\Vert \Psi_{j+1} - \Psi_j\Vert_{2}$, with $\Vert \Psi\Vert_{2}$ defined as the Frobenius norm of $ G \cdot \Psi \cdot G^*$, which represents the $L_2$-norm of $\psi$ when $\psi \in {\mathbb{H}}_{16,2}$ and $\psi|_{\mathcal{T}} = \Psi$. The rapid convergence (from $10^{-3}$ to machine precision, $10^{-14}$, in 30 seconds) indicates an efficient gradient descent optimization for the Gross-Pitaevskii ground state. The attainment of near-machine precision residuals further confirms the numerical stability of the chosen discretization scheme.

\subsubsection{Sensitivity to the initial seed II}
In the previous subsection, we applied a gradient descent method starting from an initial seed that was close to an eigenfunction of the linear part of the Hamiltonian. However, if we instead initialize with the ground state of the linear terms of the Hamiltonian, denoted as $h_{0,0}$, the optimization process will converge to the minimum more rapidly. In this example, we adopt $h_{0,0}$ as the initial seed, employing the time step $\tau=0.01$ and the same fourth-order method as before.

\begin{figure}[h!]
\centering
\begin{subfigure}[h]{0.49\linewidth}

\includegraphics[width=\linewidth]%
{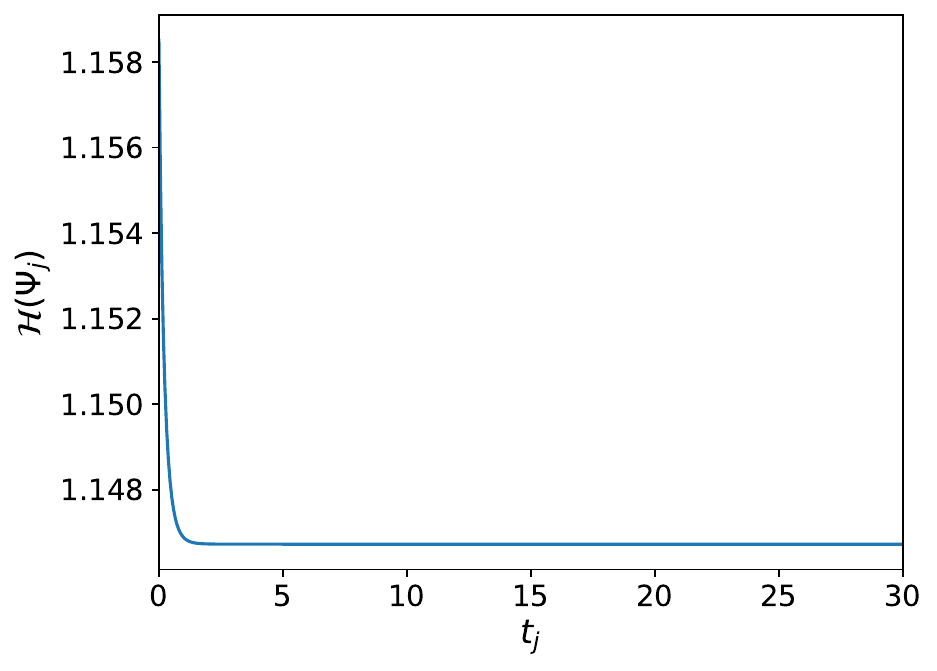}
\caption{Hamiltonian starting from the initial \\ seed $h_{0,0}$.}\label{fig_example2_Hamiltonian_faster}
\end{subfigure}
\hfill
\begin{subfigure}[h]{0.49\linewidth}
\includegraphics[width=\linewidth]%
{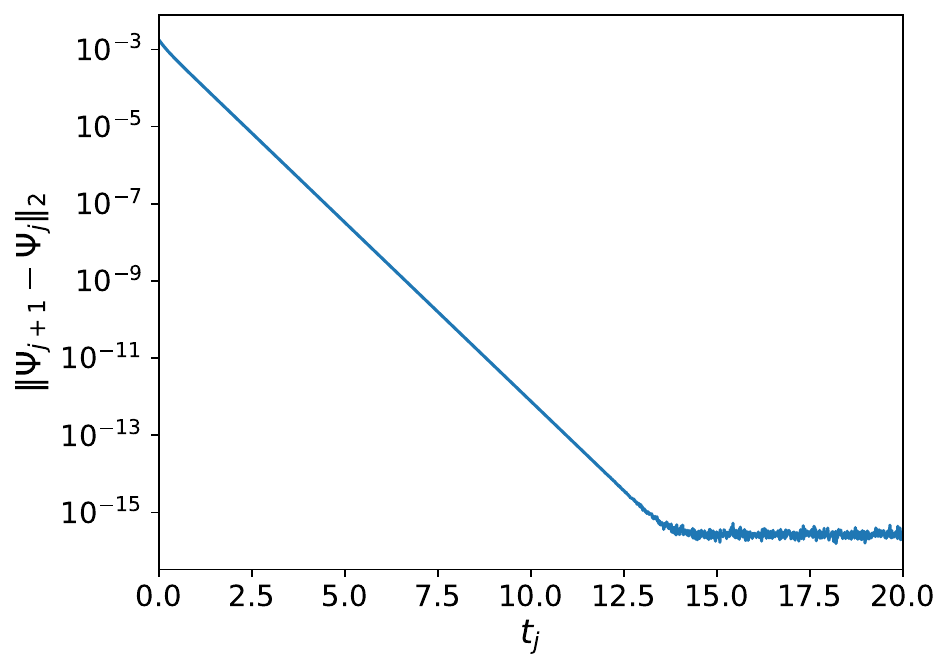}

\caption{Error between consecutive approximations of the gradient descent method.}\label{fig:example2_step_error_faster}
\end{subfigure}
\caption{Evolution of Hamiltonian and convergence of the gradient-descent iterations initialized with seed $h_{0,0}$.}
\label{fig:example2}
\end{figure}

As shown in Fig.~\ref{fig:example2}, the numerical optimization converges significantly faster compared to the previous example. The distance between successive iterations reaches machine-error precision ($\mathcal{O}(10^{-15})$) within just 15 seconds, indicating robust numerical stability. The minimum energy achieved is $\mathcal{H}_{\text{min}} = 1.146729989848536$, and the corresponding ground states agree to within $10^{-14}$, confirming numerical convergence at machine precision.

The approximations to the ground state obtained with this method are more sensitive to the chosen time step than to the order of the method. Achieving a more precise numerical ground state, approaching the exact ground state, requires progressively reducing the time step. Using this approach, the numerical exact minimum of the Hamiltonian, accurate up to $12$ decimal places, is $\mathcal{H}_{\text{min}}=1.146728491833$.

\subsection{Benchmark II: Invariant preservation under time evolution}
\label{subsec:Invariant_preservation}

We use ground states computed via the time-splitting gradient descent method for each value of the mass constraint $c$ to study the numerical behavior of the solution under the GPE time evolution \eqref{eq:Evol_Gross-Pitaevskii}. Throughout this benchmark, we evolve each ground state using a second-order splitting scheme with time step $\tau = 10^{-3}$; in the last experiment, an eigth-order scheme is also considered.

\begin{figure}[h!]

\begin{subfigure}[t]{0.48\linewidth}
\includegraphics[width=\linewidth]{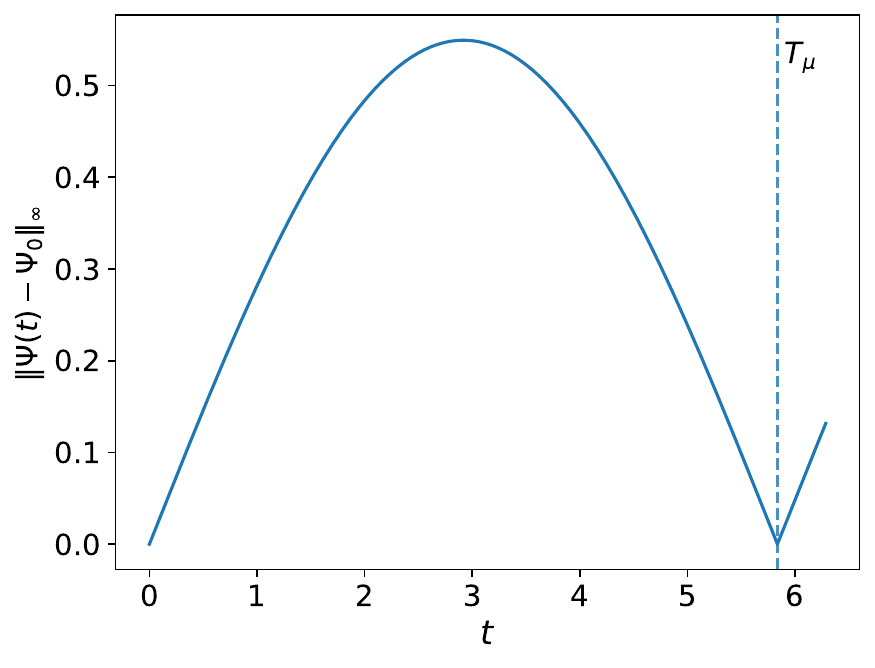}
\caption{c = 0.25, $T_\mu \approx 5.834$.}
\label{fig_Gs_vs_psi_0_c0p25}
\end{subfigure}
\hfill
\begin{subfigure}[t]{0.48\linewidth}
\includegraphics[width=\linewidth]{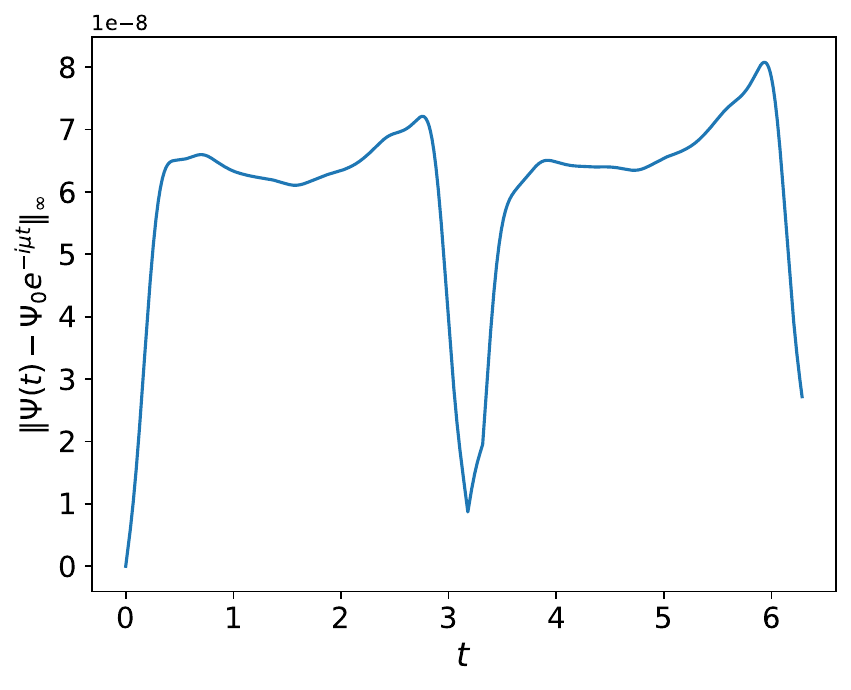}
\label{fig_Gs_vs_rotating_psi_0_norm_error_c0p25}
\vspace{-4mm}
\caption{c = 0.25.}
\end{subfigure}


\begin{subfigure}{0.48\linewidth}
\includegraphics[width=\linewidth]{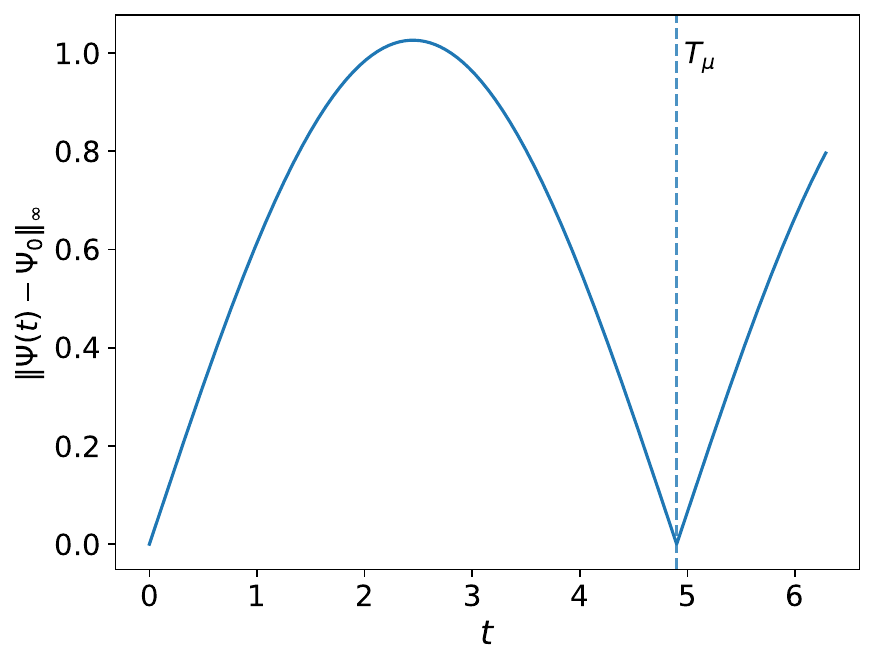}
\caption{c = 1, $T_\mu \approx 4.898$.}
\label{fig_Gs_vs_psi_0_c1p00}
\end{subfigure}
\hfill
\begin{subfigure}{0.48\linewidth}
\includegraphics[width=\linewidth]{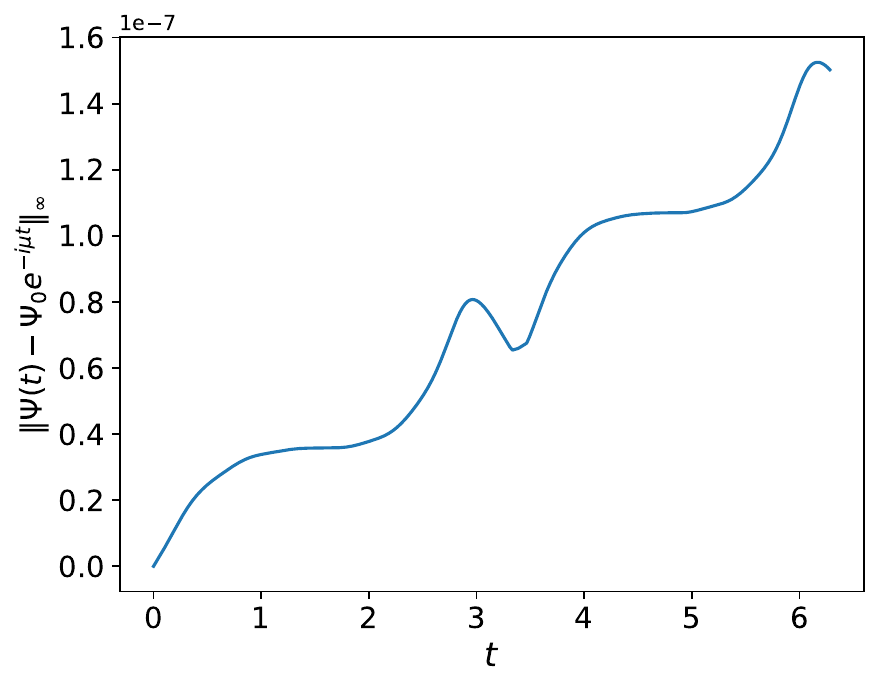}
\caption{c = 1.}
\label{fig_Gs_vs_rotating_psi_0_norm_error_c1p00}
\end{subfigure}


\begin{subfigure}{0.48\linewidth}
\includegraphics[width=\linewidth]{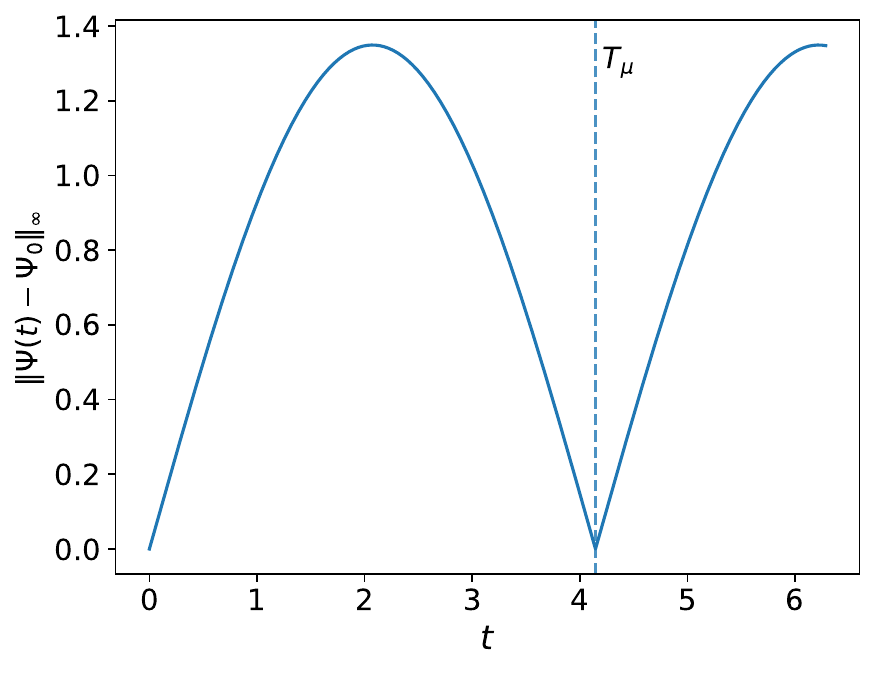}
\caption{c = 2, $T_\mu \approx 4.144$}
\label{fig_Gs_vs_psi_0_c2p00}
\end{subfigure}
\hfill
\begin{subfigure}{0.48\linewidth}
\includegraphics[width=\linewidth]{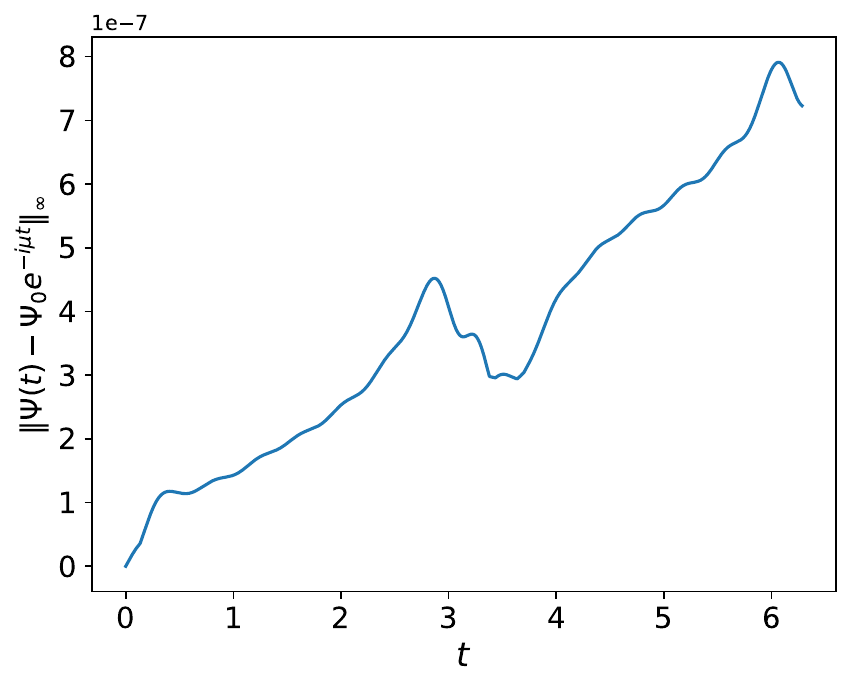}
\caption{c = 2.}
\label{fig_Gs_vs_rotating_psi_0_norm_error_c2p00}
\end{subfigure}

\caption{The panels on the left display the distance between the numerical solution $\psi(t_j)$ and the initial ground state $\psi_0$, i.e., $\|\psi(t_j) - \psi_0\|_{\infty}$, while the panels on the right show the distance to the phase-rotating ground state $\psi_0 e^{-i\mu t_j}$, i.e., $\|\psi(t_j) - \psi_0 e^{-i\mu t_j}\|_{\infty}$. Second-order splitting scheme with time step $\tau = 10^{-3}$. 
}
\label{fig_Gs_vs_psi_0}

\end{figure}

Fig.~\ref{fig_Gs_vs_psi_0} displays the distance from the evolving solution to the initial ground state, $\|\psi(t_j)-\psi_0\|_\infty$ (left panels), together with the phase-corrected distance, $\|\psi(t_j)-\psi_0 e^{-i\mu t_j}\|_\infty$ (right panels), for different mass restrictions: $c=0.25$ (first row), $c=1$ (second row) and $c=2$ (third row). We compute the theoretical period $T_{\mu}= \frac{2\pi}{\mu_{\rm ch}}$ and observe that in the three cases coincides with the numerical period of $\psi$, and the numerically computed ground state remains close to the initial ground state, undergoing a pure phase rotation.
From the left panels, one can observe that the period $T_{\mu}$ decreases as $c$ increases: $T_{\mu} = 5.834$ for $c=0.25$, $T_{\mu} = 4.898$ for $c=1$, and $T_{\mu} = 4.144$ for $c=2$. This reflects the increasing strength of the nonlinearity in the equation as the mass $c$ grows. Moreover, the right panels show a slight increase in the error as the nonlinearity becomes stronger.

\begin{figure}[t]
\centering
\begin{subfigure}[h]{0.49\linewidth}
\includegraphics[width=\linewidth]{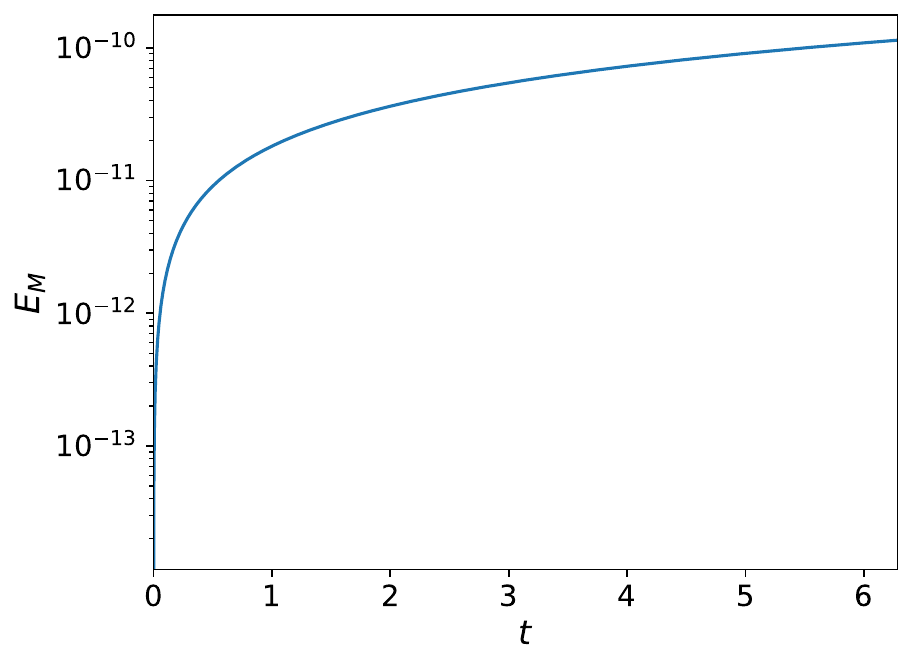}
\caption{Mass conservation: $E_M$.}\label{fig:conserv_Crel_c1_o2_h1e-3_beta2gamma1}
\end{subfigure}
\hfill
\begin{subfigure}[h]{0.49\linewidth}
\includegraphics[width=\linewidth]{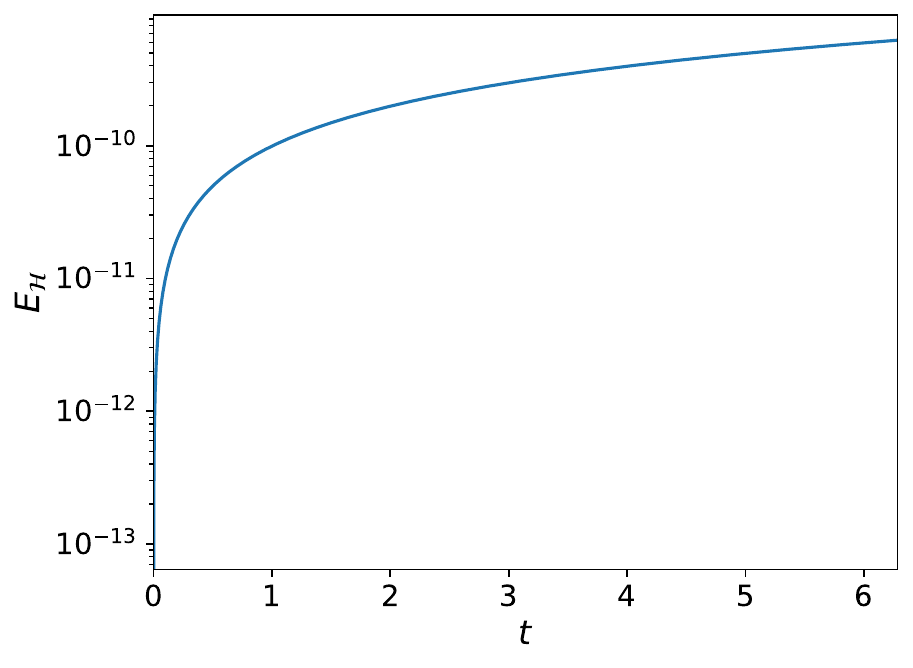}
\caption{Hamiltonian conservation: $E_\mathcal{H}$.}\label{fig:conserv_Hrel_c1_o2_h1e-3_beta2gamma1}
\end{subfigure}
\caption{Conservative properties. Method order: $q = 2$, time step: $\tau = 10^{-3}$.}
\label{fig:gs_evol_mass_and_hamiltonian_o2}
\end{figure}

In Fig.~\ref{fig:gs_evol_mass_and_hamiltonian_o2}, we monitor the conservation of mass (left) and Hamiltonian energy (right) for the ground state with mass constraint $c=1$, using a second-order time-splitting evolution method with a time step of $\tau = 10^{-3}$. We show the evolution of the relative error in mass and Hamiltonian, denoted respectively by
\[
E_M := \left\vert \frac{\left\Vert\Psi_{j}\right\Vert_2-\left\Vert\Psi_{0}\right\Vert_2}{ \left\Vert\Psi_{0}\right\Vert_2} \right\vert \quad \text{and} \quad
E_{\mathcal{H}} := \left\vert \frac{\mathcal{H}(\Psi_{j}) - \mathcal{H}(\Psi_{0})}{ \mathcal{H}(\Psi_{0})} \right\vert.
\]

Similar behaviors are observed when varying the mass constraint 
$c$, with errors remaining of the same order of magnitude. We do not report these results for brevity. In contrast, a substantial improvement in the relative errors is observed when increasing the order of the splitting methods.

\begin{figure}[t]
\centering
\begin{subfigure}[h]{0.49\linewidth}
\includegraphics[width=\linewidth]
{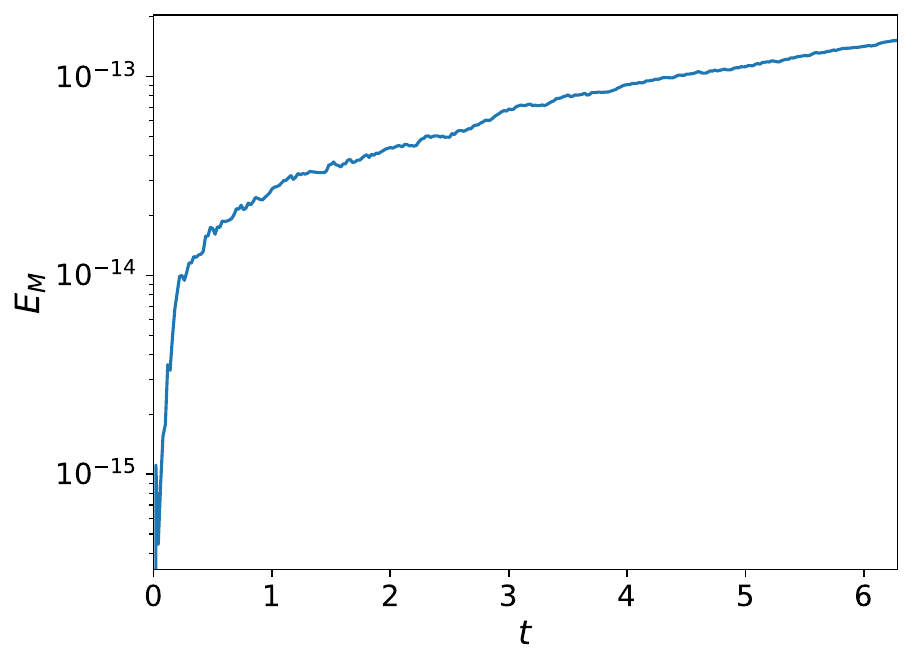}
\caption{Mass conservation: $E_M$.}\label{fig:conserv_Crel_o4_h2e-2_base17}
\end{subfigure}
\hfill
\begin{subfigure}[h]{0.49\linewidth}
\includegraphics[width=\linewidth]{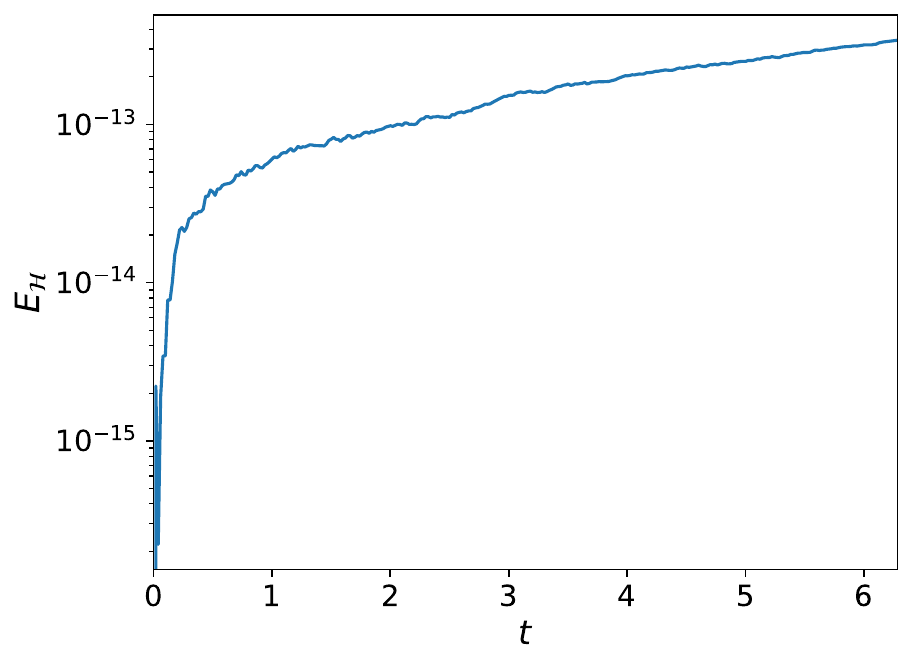}
\caption{Hamiltonian conservation: $E_{\mathcal{H}}$
.}\label{fig:conserv_Hrel_o8_h2e-2_base17}
\end{subfigure}
\caption{Conservative properties. Method order: $q = 8$, time step: $\tau = 10^{-3}$. }
\label{fig:gs_evol_mass_and_hamiltonian_o8_base17}
\end{figure}

For instance, by increasing the order of the splitting method from $q = 2$ to $q = 8$, we observe improved conservation of both mass and Hamiltonian (Fig.~\ref{fig:gs_evol_mass_and_hamiltonian_o8_base17}).

\subsection{Benchmark III: Cost–accuracy trade-off}
\label{subsec:III}
We evaluate the numerical conservation of mass and Hamiltonian energy for different method orders, for the GPE starting with the following initial condition
\[
\psi_0(x,y) = \frac{1}{2 \pi^\frac{3}{4}}e^{-V(x,y)},
\]
where $V(x,y) = \frac{x^2+y^2}{2}$ is the external trap (i.e. $\gamma=1$).
This initial condition is considered in \cite{fu2022arbitrary} because it is close to the ground state and therefore provides a convenient starting point for testing the time evolution. In our case, however, we are not restricted to this constrain, since the methodology in this work allows us to compute the ground state. Nevertheless, we use this initial condition for the sake of comparison. 

\begin{figure}[h]
\begin{center}
    \begin{subfigure}{0.48\textwidth}
        \includegraphics[width=\textwidth]{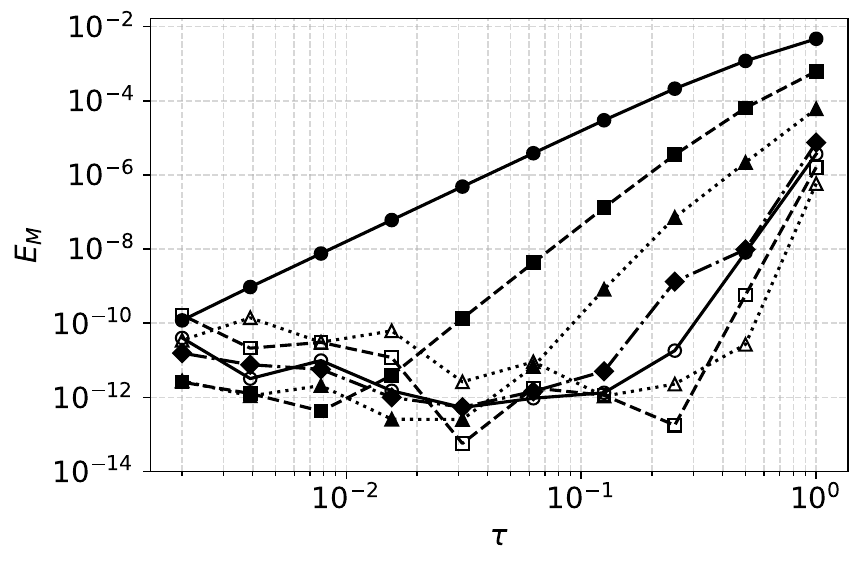}
        \caption{Relative mass error $E_M$.}
        \label{fig:cargas}
    \end{subfigure}
    \hfill
    \begin{subfigure}{0.48\textwidth}
        \includegraphics[width=\textwidth]{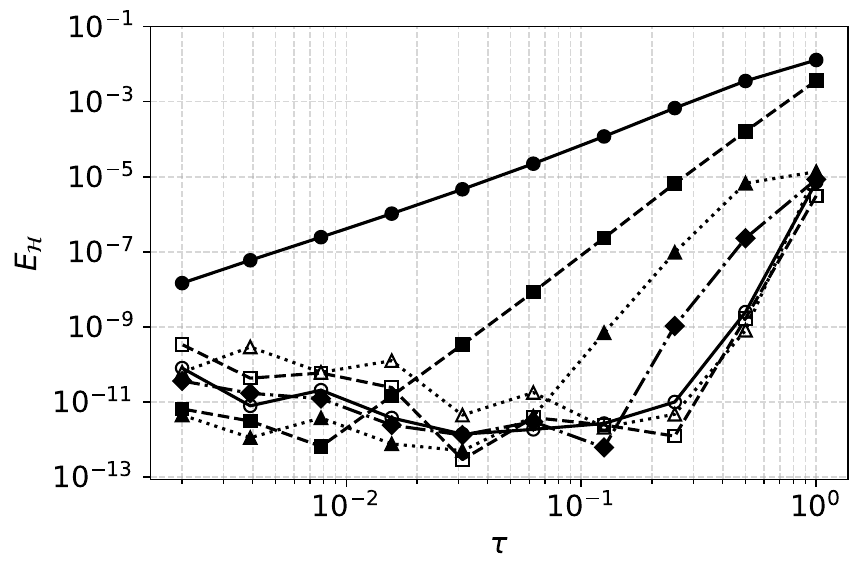}
        \caption{ Relative Hamiltonian error $E_\mathcal{H}$.}
        \label{fig:hamiltonianos}
    \end{subfigure}\\[0.5em]
    \begin{subfigure}{0.48\textwidth}
        \centering
        \includegraphics[width=\textwidth]{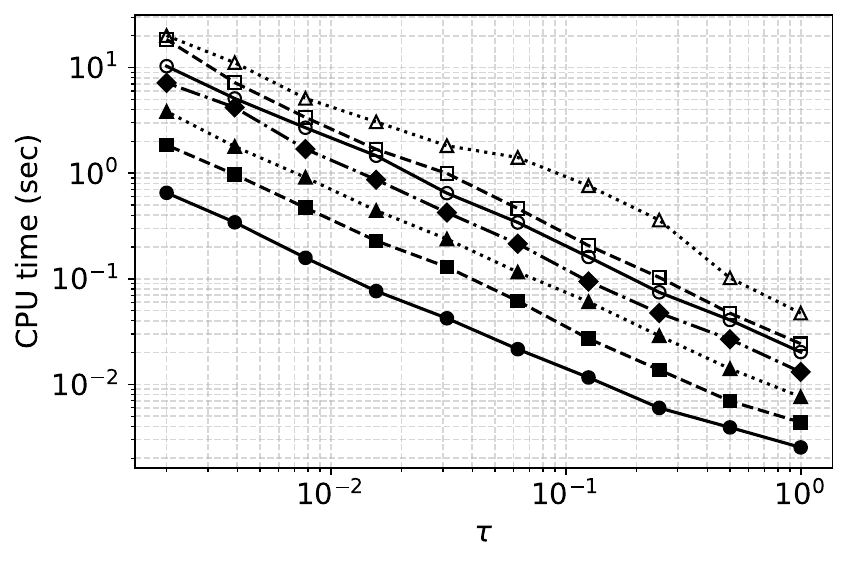}
        \caption{CPU time in seconds.}
        \label{fig:tiempos}
    \end{subfigure}
    \caption{order 2 (\text{---$\bullet$---}), order 4 (\text{- - $\scriptstyle\blacksquare$ - -}), order 6 ($\cdot \cdot \blacktriangle \cdot \cdot$), order 8 (- $\cdot {\scriptstyle\blacklozenge} \cdot$ -), order 10 (\text{---$\circ$---}), order 12 (\text{- - ${\scriptstyle\square}$ - -}), order 14 ($\cdot \cdot {\scriptstyle\triangle} \cdot \cdot$).}
    \label{fig:errores_orden}
\end{center}
\end{figure}

In Fig.~\ref{fig:errores_orden} we show the computational performance and conservation properties of the proposed splitting methods applied to the Gaussian initial condition $\psi_0$ defined above, for orders $q = 2, 4, \ldots, 14$ and final simulation time $T = 30$. Figs.~\ref{fig:cargas} and~\ref{fig:hamiltonianos} display the maximum relative errors in the mass $E_M$ and in the Hamiltonian $E_{\mathcal{H}}$, respectively, as a function of the time step $\tau$ ($\tau = \frac{1}{2^k}$, $k \in \{0, 1, \dots, 9\}$), while Fig.~\ref{fig:tiempos} shows the corresponding CPU times. All computations were performed on a Linux workstation with an AMD Ryzen 5 5600G CPU and 8 GB of RAM.

For order $q = 2$, both errors decay monotonically as $\tau$ decreases, with slopes consistent with the expected order of convergence. For higher-order methods, the errors decrease much more rapidly, reaching the round-off floor (around $10^{-12}$--$10^{-14}$) for moderate values of $\tau$, beyond which the errors stagnate and become erratic due to floating-point arithmetic limitations. In all cases, the observed slopes agree approximately with the expected order of convergence until the round-off error dominates.

In Fig.~\ref{fig:tiempos} all curves exhibit a linear behavior in the log-log scale, consistent with the expected relation $\text{CPU} \propto \tau^{-1}$: halving the time step doubles the number of steps and hence the computational cost. Since higher-order methods achieve a prescribed accuracy with a much larger $\tau$, they can be significantly more efficient than lower-order methods when high precision is required.

In Section \ref{sec:High-order_TSM}, we already mentioned that the flows of these methods are amenable to parallelization, which would significantly reduce the computational cost, particularly for high-order methods. Since no parallel implementation is considered here, the curves exhibit a vertical displacement that grows with the order $q$ as $q\left(\frac{q}{2}+1\right)$, consistent with the theoretical analysis in Section \ref{sec:High-order_TSM} and \cite{DeLeo2015} (see Eq.~\eqref{eq:steps_order}). To illustrate this, Table~\ref{table:tiempos_orden} shows the ratio of the mean CPU time for each order with respect to order $q=2$. The observed values agree well with the theoretical prediction for low orders, while some deviation is expected for higher orders due to computational overhead not accounted for in the theoretical model.

The observed ranges for the relative errors in mass and energy conservation, as well as the computational times, are consistent with those reported in the literature (cf. \cite{fu2022arbitrary}), where methods up to order 6 are considered for the case $V = 0$. In contrast, the present work extends the numerical assessment to orders up to $q = 14$ in the presence of a harmonic trapping potential, suggesting that the proposed methods offer competitive performance even under limited computational resources, as all results here were obtained on a standard desktop workstation without parallelization or specialized hardware.

\begin{table}[h]
\centering
\begin{tabular}{ccc}
\hline
Order $q$ & Observed ratio & Theoretical ratio \\
\hline
2  & 1.00  & 1.00  \\
4  & 2.99  & 3.00  \\
6  & 5.80  & 6.00  \\
8  & 11.34 & 10.00 \\
10 & 19.15 & 15.00 \\
12 & 22.02 & 21.00 \\
14 & 40.07 & 28.00 \\
\hline
\end{tabular}
\caption{Ratio of mean CPU time with respect to order 2, observed and theoretical.}
\label{table:tiempos_orden}
\end{table}

\section{Discussion and Conclusions}

In this work, we presented a high-order time-splitting framework for the numerical approximation of the two-dimensional GPE, addressing both the computation of the ground state and the time evolution of the solution. The methods were designed to ensure conservation of key physical quantities, such as mass and Hamiltonian energy, while maintaining computational efficiency.

We demonstrated that combining gradient descent with time-splitting methods provides an effective strategy for obtaining accurate approximations of the ground state. The convergence behavior, analyzed for different initial seeds, highlights the sensitivity of the method to the initial condition and time step size. In particular, initializing the method with the linear ground state significantly accelerates convergence and improves stability.

The time evolution tests further confirmed the robustness of the proposed methods. Across a wide range of method orders and time steps, the relative errors in mass and energy conservation remained within acceptable bounds. The numerical experiments in Section \ref{subsec:III} validate that higher-order methods achieve accurate results even with larger time steps, supporting their use in computationally demanding simulations.

It is important to note that all computations were carried out on a standard notebook computer, without the use of parallelization or specialized hardware. Therefore, the reported computational times are provided for reference purposes only. Despite this, the observed ranges for the relative errors in mass and energy conservation, as well as the computational times, are consistent with those reported in the literature (cf., \cite{fu2022arbitrary}). This agreement reinforces the validity and competitiveness of the proposed approach, even under modest computational resources.

Overall, the numerical evidence supports the use of high-order positive time-splitting methods for efficiently and accurately solving the GPE. Future work may explore adaptive strategies that balance splitting order, basis size, and time step to further optimize accuracy and cost, as well as extensions to more complex or higher-dimensional configurations.

\section*{Acknowledgments}
This work was supported by the following projects: PI UNGS 30/3394 and CONICET PIP 11220220100124CO.
\bigskip


\section*{Declaration of generative AI and AI-assisted technologies in the manuscript preparation process}

During the preparation of this work the authors used Claude (Anthropic) in order to assist in the revision and editing of the written presentation of the manuscript. After using this tool, the authors reviewed and edited the content as needed and take full responsibility for the content of the published article.

\bibliography{biblio2}

@article{antoine2013computational,
  title={Computational methods for the dynamics of the nonlinear {S}chr{\"o}dinger {G}ross--{P}itaevskii equations},
  author={Antoine, Xavier and Bao, Weizhu and Besse, Christophe},
  journal={Computer Physics Communications},
  volume={184},
  number={12},
  pages={2621--2633},
  year={2013},
  publisher={Elsevier}
}

@misc{assanto2012nematicons,
  title={Nematicons, spatial optical solitons in nematic liquid crystals},
  author={Assanto, G},
  year={2012},
  publisher={Wiley, New York}
}

@article{bao2004computing,
  title={Computing the ground state solution of {B}ose--{E}instein condensates by a normalized gradient flow},
  author={Bao, Weizhu and Du, Qiang},
  journal={SIAM Journal on Scientific Computing},
  volume={25},
  number={5},
  pages={1674--1697},
  year={2004},
  publisher={SIAM}
}

@article{BaoCai2013,
  author  = {Bao, Weizhu and Cai, Yongyong},
  title   = {Mathematical theory and numerical methods for {B}ose--{E}instein condensation},
  journal = {Kinetic and Related Models},
  volume  = {6},
  number  = {1},
  pages   = {1--135},
  year    = {2013},
  doi     = {10.3934/krm.2013.6.1}
}

@article{bao2009generalized,
  title={A generalized-{L}aguerre--{F}ourier--{H}ermite pseudospectral method for computing the dynamics of rotating {B}ose--{E}instein condensates},
  author={Bao, Weizhu and Li, Hailiang and Shen, Jie},
  journal={SIAM Journal on Scientific Computing},
  volume={31},
  number={5},
  pages={3685--3711},
  year={2009},
  publisher={SIAM}
}

@article{ben2015localized,
  title={Localized solutions for a nonlocal discrete {NLS} equation},
  author={Ben, Roberto I and Ake, Luis Cisneros and Minzoni, AA and Panayotaros, Panayotis},
  journal={Physics Letters A},
  volume={379},
  number={30-31},
  pages={1705--1714},
  year={2015},
  publisher={Elsevier}
}

@article{ben2017properties,
  title={Properties of some breather solutions of a nonlocal discrete {NLS} equation},
  author={Ben, Roberto Ignacio and Borgna, Juan Pablo and Panayotaros, Panayotis},
  journal={Communications in Mathematical Sciences},
volume={15},
number={8},
pages={2143-2175},
  year={2017},
  publisher={International Press Boston}
}

@article{besse2002order,
  title={Order estimates in time of splitting methods for the nonlinear {S}chr{\"o}dinger equation},
  author={Besse, Christophe and Bid{\'e}garay, Brigitte and Descombes, St{\'e}phane},
  journal={SIAM Journal on Numerical Analysis},
  volume={40},
  number={1},
  pages={26--40},
  year={2002},
  publisher={SIAM}
}

@article{besteiro1,
  title={Existence of Peregrine type solutions in fractional reaction–diffusion equations.},
  author={A. Besteiro and D. Rial},
  journal={Electron. J. Qual. Theory Differ. Equ.},
  	Pages = {1-9},
  Volume = {9},
	Year = {2019}}

@article{Besteiro2018,
  title={Global existence for vector valued fractional reaction-diffusion equations},
  author={A. Besteiro and D. Rial},
  journal={Publicacions Matemàtiques},
  Pages = {653-680},
  Volume = {96},
  Issue = {2},
  year={2021}
}

@article{Borgna2015,
	Author = {J. P. Borgna and M. De Leo and D. Rial and C. S\'anchez de la Vega},
	Journal = {Commun. Math. Sci, Int. Press Boston, Inc.},
	Pages = {83-101},
	Title = {General Splitting methods for abstract semilinear evolution equations},
	Volume = {13},
	Year = {2015}}

@article{borgna2018optical,
  title={Optical solitons in nematic liquid crystals: model with saturation effects},
  author={Borgna, Juan Pablo and Panayotaros, Panayotis and Rial, Diego and De la Vega, Constanza S{\'a}nchez F},
  journal={Nonlinearity},
  volume={31},
  number={4},
  pages={1535},
  year={2018},
  publisher={IOP Publishing}
}

@book{carretero2024nonlinear,
  title={Nonlinear Waves and Hamiltonian Systems: From One to Many Degrees of Freedom, from Discrete to Continuum},
  author={Carretero-Gonz{\'a}lez, Ricardo and Frantzeskakis, Dimitrios J and Kevrekidis, Panayotis G},
  year={2024},
  publisher={Oxford University Press}
}

@article{castella2009splitting,
  title={Splitting methods with complex times for parabolic equations},
  author={Castella, Fran{\c{c}}ois and Chartier, Philippe and Descombes, St{\'e}phane and Vilmart, Gilles},
  journal={BIT Numerical Mathematics},
  volume={49},
  pages={487--508},
  year={2009},
  publisher={Springer}
}

@article{DeLeo2015,
	Author = {M. {De Leo} and D. Rial and C. F. S\'anchez de la Vega},
	Journal = {IMA J. Numer. Anal.},
        volume = {36},
        number = {4},
	Pages = {1842-1866},
	Title = {High-order time-splitting methods for irreversible equations},
	Year = {2016}}

@article{descombes2010exact,
  title={An exact local error representation of exponential operator splitting methods for evolutionary problems and applications to linear {S}chr{\"o}dinger equations in the semi-classical regime},
  author={Descombes, St{\'e}phane and Thalhammer, Mechthild},
  journal={BIT Numerical Mathematics},
  volume={50},
  number={4},
  pages={729--749},
  year={2010},
  publisher={Springer}
}

@article{descombes2013lie,
  title={The {L}ie--{T}rotter splitting for nonlinear evolutionary problems with critical parameters: a compact local error representation and application to nonlinear {S}chr{\"o}dinger equations in the semiclassical regime},
  author={Descombes, St{\'e}phane and Thalhammer, Mechthild},
  journal={IMA Journal of Numerical Analysis},
  volume={33},
  number={2},
  pages={722--745},
  year={2013},
  publisher={Oxford University Press}
}

@article{fu2022arbitrary,
  title={Arbitrary high-order exponential integrators conservative schemes for the nonlinear {G}ross-{P}itaevskii equation},
  author={Fu, Yayun and Hu, Dongdong and Zhang, Gengen},
  journal={Computers \& Mathematics with Applications},
  volume={121},
  pages={102--114},
  year={2022},
  publisher={Elsevier}
}

@article{gauckler2011convergence,
  title={Convergence of a split-step {H}ermite method for the {G}ross--{P}itaevskii equation},
  author={Gauckler, Ludwig},
  journal={IMA journal of Numerical Analysis},
  volume={31},
  number={2},
  pages={396--415},
  year={2011},
  publisher={OUP}
}

@article{gross1961structure,
  title={Structure of a quantized vortex in boson systems},
  author={Gross, E. P.},
  journal={Il Nuovo Cimento},
  volume={20},
  pages={454--477},
  year={1961}
}

@book{kevrekidis2008emergent,
  title={Emergent nonlinear phenomena in Bose-Einstein condensates: theory and experiment},
  author={Kevrekidis, Panayotis G and Frantzeskakis, Dimitri J and Carretero-Gonza{\'a}lez, R},
  volume={45},
  year={2008},
  publisher={Springer}
}

@book{kevrekidis2015defocusing,
  title={The defocusing nonlinear Schr{\"o}dinger equation: from dark solitons to vortices and vortex rings},
  author={Kevrekidis, Panayotis G and Frantzeskakis, Dimitri J and Carretero-Gonz{\'a}lez, Ricardo},
  year={2015},
  publisher={SIAM, Society for Industrial and Applied Mathematics}
}

@article{lubich2008splitting,
  title={On splitting methods for {S}chr{\"o}dinger-{P}oisson and cubic nonlinear {S}chr{\"o}dinger equations},
  author={Lubich, Christian},
  journal={Mathematics of computation},
  volume={77},
  number={264},
  pages={2141--2153},
  year={2008}
}

@article{lieb2001rigorous,
  title={A Rigorous Derivation of the {G}ross--{P}itaevskii {E}nergy {F}unctional for a {T}wo-dimensional {B}ose {G}as},
  author={Lieb, Elliott H and Seiringer, Robert and Yngvason, Jakob},
  journal={Communications in Mathematical Physics},
  volume={224},
  pages={17--31},
  year={2001},
  publisher={Springer}
}

@article{neri1987lie,
  title={Lie algebras and canonical integration},
  author={Neri, Filippo},
  journal={Dept. of Physics, University of Maryland},
  year={1987}
}

@book{pelinovsky2011localization,
  title={Localization in periodic potentials: from {S}chr{\"o}dinger operators to the {G}ross--{P}itaevskii equation},
  author={Pelinovsky, Dmitry E},
  volume={390},
  year={2011},
  publisher={Cambridge University Press}
}

@article{pitaevskii1961vortex,
  title={Vortex lines in an imperfect Bose gas},
  author={Pitaevskii, Lev P},
  journal={Sov. Phys. JETP},
  volume={13},
  number={2},
  pages={451--454},
  year={1961}
}

@book{pethick2008bose,
  title={Bose--Einstein condensation in dilute gases},
  author={Pethick, Christopher J and Smith, Henrik},
  year={2008},
  publisher={Cambridge university press}
}

@article{ruth1983canonical,
  title={A canonical integration technique},
  author={Ruth, Ronald D},
  journal={IEEE Trans. Nucl. Sci.},
  volume={30},
  number={CERN-LEP-TH-83-14},
  pages={2669--2671},
  year={1983}
}

@article{yoshida1990construction,
  title={Construction of higher order symplectic integrators},
  author={Yoshida, Haruo},
  journal={Physics letters A},
  volume={150},
  number={5-7},
  pages={262--268},
  year={1990},
  publisher={Elsevier}
}
\bibliographystyle{siam}

\end{document}